\documentclass[11pt, a4paper, fleqn, twoside]{article}
\usepackage{amsmath}
\usepackage{amsfonts}
\usepackage{eufrak}
\usepackage{amsthm}
\usepackage{amssymb}
\usepackage[all]{xy}
\SelectTips{cm}{}

\title{On the algebraic classification of $K$-local spectra}
\author{Constanze Roitzheim}
\date{}

\DeclareMathOperator{\Hom}{Hom}
\DeclareMathOperator{\Ho}{Ho}
\DeclareMathOperator{\Ext}{Ext}
\DeclareMathOperator{\ma}{max}
\DeclareMathOperator{\injdim}{injdim}
\DeclareMathOperator{\md}{mod}
\DeclareMathOperator{\ke}{ker}
\DeclareMathOperator{\im}{im}
\DeclareMathOperator{\coker}{coker}
\DeclareMathOperator{\el}{else}
\DeclareMathOperator{\colim}{colim}
\DeclareMathOperator{\odd}{odd}
\DeclareMathOperator{\even}{even}
\DeclareMathOperator{\Hocolim}{Hocolim}
\DeclareMathOperator{\cone}{cone}
\DeclareMathOperator{\map}{map}

\DeclareMathOperator{\sset}{sSet_*}
\DeclareMathOperator{\fib}{fib}

\DeclareMathOperator{\fatl}{\bf{l}}
\DeclareMathOperator{\fatr}{\bf{r}}

\theoremstyle{definition}
\newtheorem{ddd}{Definition}[subsection]
\newtheorem{ppp}[ddd]{Proposition}
\newtheorem{lem}[ddd]{Lemma}
\newtheorem{ttt}[ddd]{Theorem}
\newtheorem{ccc}[ddd]{Corollary}
\newtheorem*{rrr}{Remark}
\newtheorem*{nnn}{Notation}

\newcommand{\vsp}{\vspace{0.7cm}}

\newcommand{\A}{\mathcal{A}}
\newcommand{\B}{\mathcal{B}}
\newcommand{\C}{\mathcal{C}}
\newcommand{\D}{\mathcal{D}}
\newcommand{\Q}{\mathcal{Q}}
\newcommand{\R}{\mathcal{R}}
\newcommand{\E}{L_1\mathcal{S}}
\newcommand{\Ca}{\mathcal{C}^1(\mathcal{A})}
\newcommand{\Cb}{\mathcal{C}^{2p-2}(\mathcal{B})}
\newcommand{\Da}{\mathcal{D}^1(\mathcal{A})}
\newcommand{\Db}{\mathcal{D}^{2p-2}(\mathcal{B})}
\newcommand{\xtild}{\tilde{X}}

\begin{document}
\maketitle

\begin{abstract}
In 1996, Jens Franke proved the equivalence of certain triangulated
categories possessing an Adams spectral sequence. One particular
application of this theorem is that the $K_{(p)}$-local stable
homotopy category at an odd prime can be described as the derived
category of an abelian category. We explain this proof from a
topologist's point of view.
\end{abstract}

\vsp
 In 1983 Bousfield published a paper about the category of
$E(1)$-local (or, equivalently, $K$-local) spectra at an odd prime.
There, he gave an algebraic description of isomorphism classes of
$E(1)$-local spectra in their homotopy category via $E(1)$-homology
and a certain ``$k$-invariant'' coming from a $d_2$-differential in
the Adams spectral sequence. However, with this setup he could only
describe the morphisms up to Adams filtration.

In 1996, Jens Franke constructed an abstract equivalence between
certain triangulated categories possessing an Adams spectral
sequence. Applying Franke's main theorem to the special case of
$E(1)$-local spectra, one obtains an algebraic description of the
homotopy category of $E(1)$-local spectra also covering the
morphisms. In this paper, we give a streamlined exposition of
Franke's result adapted to this special case:

\vsp
{\bf Theorem}[Franke] There is an equivalence of categories
\[
\R: \Db \longrightarrow \Ho(\E)
\]
where $\Db$ denotes the derived category of twisted cochain complexes over the abelian category $\B$, and
$\Ho(\E)$ the homotopy category of $E(1)$-local spectra.
\vsp

This paper is organised as follows: In the first chapter, the
categories playing the main role for the construction are
introduced: firstly, the category of so-called twisted cochain
complexes of $E(1)_* E(1)$- comodules and secondly, a certain
diagram category of spectra with a fixed diagram shape and a model
structure related to the model structure of $E(1)$-local spectra.

In the next section, a functor $\Q$ is constructed which gives an equivalence of twisted
cochain complexes and the homotopy category of above diagram spectra. In the third section
this equivalence $\Q$ is extended to an equivalence of the derived category of twisted cochain
complexes and the homotopy category of $E(1)$-local spectra. Further, as section 4 will show,
this equivalences gives an ``exotic model'' for $E(1)$-local spectra: the homotopy
categories of the cochain complexes and $E(1)$-local spectra are equivalent as categories,
yet there is no Quillen equivalence between them.

We do not claim any originality, it is just the proof of Franke's Main Uniqueness
Theorem applied to Bousfield's case with the notation adapted and some technical details
filled in. My special thanks go to Stefan Schwede for his motivation and support.

\vsp
\section{The main ingredients}

\subsection{$E(1)_\ast E(1)$-comodules}

We begin with describing an abelian category denoted $\A$ which is equivalent to the category of
$E(1)_{\ast} E(1)$-comodules (see \cite{Bou85}, 10.3).
Bousfield describes $\A$ as follows: Let $p$ be an odd prime and let $\B = \B(p)$ denote the
category of $\mathbb{Z}_{(p)}$-modules together with Adams operations $\psi^k, \quad
k \in \mathbb{Z}_{(p)}^\ast$ satisfying the following:

For each $M \in \B(p)$,
\begin{itemize}

\item There is an eigenspace decomposition

\[
M \otimes \mathbb{Q} \cong \bigoplus\limits_{j \in \mathbb{Z}} W_{j(p-1)}
\]

such that for all $w \in W_{j(p-1)}$ and $k \in \mathbb{Z}_{(p)}$:

\[
(\psi^k \otimes id)w = k^{j(p-1)}w.
\]

\item For all $x \in M$ there is a finitely generated submodule $C(x)$ containing $x$ satisfying:
 for all $m \ge 1$ there is an $n$ such that the action of $\mathbb{Z}_{(p)}^\ast$ on
$C(x) \big/ p^m C(x)$ factors through the quotient of $({ \mathbb{Z}} \big/
{\scriptstyle (p^{n+1})})^\ast$ by
its subgroup of order $p-1$.

\end{itemize}

To build the category $\A$ out of the above category, we additionally need the following:
Let $T^{j(p-1)}: \B \longrightarrow \B, j\in\mathbb{Z}$, denote the following self-equivalence:

For all $M \in \B, \quad T^{j(p-1)}(M) = M \quad$ as a $\mathbb{Z}_{(p)}$-module, but on $T^{j(p-1)}(M)$,
the Adams operation $\psi^k$ now equals $k^{j(p-1)}\psi^k: M \longrightarrow M$ for all
$k\in\mathbb{Z}$.

Now an object $\mathcal{M}\in\A$ is defined as a collection of modules
$\mathcal{M} = (M_n)_{n\in\mathbb{Z}}, M_n\in\B$, together with isomorphisms
\[
T^{p-1}(M_n) \longrightarrow M_{n+2p-2} \quad\mbox{for all}\quad n \in \mathbb{Z}.
\]

In this paper we will often make use of the following:
Let $X$ be a spectrum. Then the $E(1)_*E(1)$-comodule $E(1)_*(X)$ is an object of $\A$ in the above sense
by taking $M_n:=E(1)_n(X)$, and the operations $\psi^k$ being the usual Adams operations.

From now on $\B$ will be viewed as the subcategory of $\A$ consisting of those objects
$(M_n)_{n\in\mathbb{Z}}$ such that
\[
M_n = \left\{ \begin{array}{r@{\quad:\quad}l}
                 M & n \equiv 0 \quad \md  \quad 2p-2 \\ 0 & \el
                \end{array} \right.
\]

This describes a so-called {\it split} of period $2p-2$ of $\A$: $\B
\subset \A$ is a Serre class such that
\[
\hspace{3cm}
\bigoplus\limits_{0 \le i < 2p-2} \B \longrightarrow \A
\]
\[
\hspace{2.7cm}
(B_i)_{0 \le i < 2p-2} \longmapsto \bigoplus\limits_{0 \le i < 2p-2}B_i[i]
\]
is an equivalence of categories, where $[i]$ denotes the $i$-fold internal shift in the
grading, i.e.
$ M[i]_n = M_{n-i}. $

\begin{rrr}
There exists a similar splitting of period $2p-2$ for the category of
$E(n)_\ast E(n)$-comodules with arbitrary $n$ and $p$ odd. Moreover, the proof of the uniqueness theorem
will not only work for the case $p$ odd and $n=1$ but for all $p$ and $n$ such that
$n^2+n < 2p-2$, i.e. when the maximal injective dimension of $E(n)_\ast E(n)$-comodules
is smaller than the splitting period.
\end{rrr}

\subsection{Twisted cochain complexes}

In this section we describe the source of the equivalence to be constructed. Let $\A$ for
the next paragraphs denote an arbitrary abelian category, $N$ a natural number and \linebreak
$ T: \A \longrightarrow \A$ a self-equivalence.

\begin{ddd} The category $\C^{(T,N)}(\A)$ of $(T,N)$-{\it twisted cochain complexes} with \linebreak values in
$\A$ is defined as follows:

The objects are cochain complexes $C^\ast$ with $C^i \in \A$ for all $i$ together with
an isomorphism of cochain complexes
\[
\alpha_C: T(C^\ast) \longrightarrow C^\ast[N]=C^{\ast+N}.
\]
The morphisms are morphisms of cochain complexes $f: C^\ast \rightarrow D^\ast$ that are
compatible with those isomorphims, i.e. the following diagram commutes:
\[
\xymatrix{
T(C^\ast)\ar[r]^{\alpha_C}\ar[d]_{T(f)} & C^\ast[N]\ar[d]^{f[N]} \\
T(D^\ast)\ar[r]^{\alpha_D} & D^\ast[N].
}
\]

Such a cochain complex $C^\ast$ is called {\it injective} if each $C^i$ is injective in $\A$. A
morphism in $\C^{(T,N)}(\A)$ is called a quasi-isomorphism if it induces an isomorphims in
cohomology. $C^\ast$ is called {\it strictly injective} if it is injective and for each acyclic
complex $D^\ast$, the cochain complex $\Hom^\ast_{\C^{(T,N)}(\A)}(D^\ast,C^\ast)$ is again acyclic.
\end{ddd}

\begin{nnn}
In our particular case, let $\A$ be again the category equivalent to
$E(1)_\ast E(1)$-comodules described in the last section. The self-equivalence of $\A$ we work with
from now on is last section's $T^{p-1}.$ We denote the category $\C^{(T^{p-1},1)}(\A)$
by $\Ca$.
\end{nnn}

Secondly, we are interested in the category $\C^{(T^{(2p-2)(p-1)},2p-2)}(\B)$, where $\B$
denotes again the split of $\A$ introduced in the last section. This category of cochain
complexes will be denoted by $\Cb$.

\subsection{A model structure for twisted cochain complexes}

\vsp
\begin{ppp} {[Franke]} There is a model structure on $\Ca$ resp. $\Cb$ such that
\begin{itemize}
\item weak equivalences are the quasi-isomorphisms
\item cofibrations are the monomorphisms
\item fibrations are the degreewise split epimorphisms with strictly injective kernel.
\end{itemize}
\end{ppp}

\begin{rrr}
The analogous model structure exists on arbitrary $\C^{(T,N)}(\A)$, given that there are
enough injectives in $\A$.
\end{rrr}

\begin{nnn}
$\Da$ resp. $\Db$ denotes the derived category of $\Ca$ resp. $\Cb$, i.e.
the homotopy category of these model categories with respect to the above model structure.
\end{nnn}

\vsp
\subsection{\bf The relation between $\Ca$ and $\Cb$}

\vsp
Now we will describe in which ways $\Ca$ and $\Cb$ contain the same data and therefore
are equivalent as categories.

Let $C^\ast =( \hspace{0.5cm}...\rightarrow C^0 \rightarrow C^1 \rightarrow C^2 \rightarrow ...
\hspace{0.5cm} )$
be an object of $\Ca$, i.e. $C^i \in \A$ and $T^{p-1}(C^i)\cong C^{i+1}$ via $\alpha_C$.
Since $\A$ splits into $2p-2$ copies of $\B$, each $C^i$ splits into
$C^i=C^i_{(0)} \oplus C^i_{(1)}
\oplus ... \oplus C^i_{(2p-1)}$ with $C^i_{(j)} \in \B[j]$. So $C^*$ gives us a complex
taking values in $\B$ by setting
\[
C^\ast_{(0)} :=( \hspace{0.5cm}...\rightarrow C^0_{(0)}
\rightarrow C^1_{(0)} \rightarrow C^2_{(0)} \rightarrow ... \hspace{0.5cm} ).
\]
The self-equivalence $T^{p-1}$ acts on each $C^i$ by cyclically permuting the summands:
\[
T(C^i_{(j)}) \cong T(C^i)_{(j+1)} \cong C^{i+1}_{(j+1)}, j \in \mathbb{Z} \big/ {\scriptstyle (2p-2)}.
\]
Consequently we have
\[
T^{(2p-2)(p-1)}(C^i_{(0)}) \cong T^{(2p-3)(p-1)}(C^{i+1}_{(1)}) \cong ... \cong
C^{i+2p-2}_{(0)},
\]
and thus $C^*_{(0)}$ is $2p-2$-twistperiodic, i.e. $C^*_{(0)} \in Obj(\Cb)$.

On the other hand, an object of $\Cb$ carries the same information as an object of $\Ca$:
given
\[
D^\ast =( \hspace{0.5cm}...\rightarrow D^0 \rightarrow D^1 \rightarrow D^2 \rightarrow ...
\hspace{0.5cm} ) \in \Cb
\]
 one obtains a corresponding complex $\overline{D}^* \in \Ca$
by setting
\[
\overline{D}^i_{(j)}:=T^{j(p-1)}(D^{i-j}).
\]
 So all in all, it is of no significant
relevance which of those two categories we choose to work in.

\subsection{Diagram categories of spectra}

By a {\it spectrum} we mean the following: A spectrum $X$ is a
collection of simplicial sets $X_n$ for $n \ge 0$ together with
morphisms of simplicial sets $\sigma_n: \Sigma X_n \rightarrow
X_{n+1}$. A morphism $f: X \rightarrow Y$ of spectra is a collection
of morphisms $f_n: X_n \rightarrow Y_n$ of simplicial sets that
commute with the structure maps $\sigma_n$, i.e. $\sigma_n \circ
\Sigma f_n  = f_{n+1} \circ \sigma_n$ (see \cite{BouFri78}). Let
$\E$ denote the category of spectra together with the following
model structure which is a localisation of the Bousfield-Friedlander
model structure: $f:X \longrightarrow Y$ is a

\begin{itemize}
\item weak equivalence if $E(1)_*(f)$ is an isomorphism in $\A$
\item cofibration if each $g_n: X_n \bigcup\limits_{\Sigma X_{n-1}} \Sigma Y_{n-1}
\longrightarrow Y_n$ is a cofibration of simplicial sets.
\item fibration if $f$ has the right lifting property with respect to acyclic cofibrations
\end{itemize}

(This model structure is rather well-known, however, we do not know any reference in literature.)
Note that $\Ho(\E)$ is equivalent to the homotopy category of $E(1)$-local spectra denoted
$\Ho(L_1 \mathcal{S})$.

By a {\it poset} we mean a partially ordered finite set. For a poset $C$, $\E^C$ denotes the
category of $C$-shaped diagrams with values in $\E$. For each $c \in C$ and $X \in \E^C$,
let $X_c$ denote the value of $X$ at the vertex $c$. For example, taking the poset
$\underline{1}=(0 \rightarrow 1)$, an object of $\E^{\underline{1}}$ is determined by a morphism
$X_0 \longrightarrow X_1$ in $\E$.

For fixed $C$, there is a model structure on $\E^C$: A morphism $f: X \longrightarrow Y$
of diagrams is a

\begin{itemize}
\item a weak equivalence if it is a vertexwise weak equivalence in $\E$ (i.e. $f_c: X_c \rightarrow Y_c$ induces an isomorphism
in $E(1)$-homology for each $c \in C$)
\item a fibration if it is vertexwise a fibration in $\E$
\item a cofibration if for all $c \in C$, $X_c \coprod\limits_{\colim_{c'<c}X_{c'}} Y_c'
\longrightarrow Y_c$ is a cofibration.

\end{itemize}

This gives $\E^C$ the structure of a stable model category, thus $\Ho(\E^C)$ is a triangulated
category (see e.g. \cite{Hov99}).

\vsp
From now on, $C$ will be the poset consisting of elements $\beta_i$ and $\gamma_i$ for $i \in \mathbb{Z} \big/
{\scriptstyle (2p-2)} $ such that $\beta_i > \gamma_i$ and $\beta_i > \gamma_{i+1}$ for $i \in \mathbb{Z} \big/ {\scriptstyle (2p-2)}$, i.e.

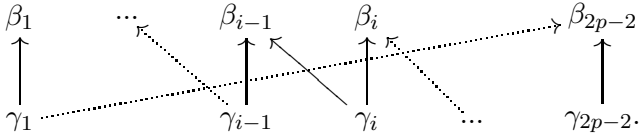
\[
\xymatrix{      \beta_1 & ... & \beta_{i-1} & \beta_i & & \beta_{2p-2} \\
             \gamma_1 \ar[u]\ar@{.>}[urrrrr] & & \gamma_{i-1}\ar[u]\ar@{.>}[ul] & \gamma_i\ar[u]\ar[ul] & ... \ar@{.>}[ul]
            & \gamma_{2p-2}. \ar[u]
}
\]

So an object $X$ of $\Ho(\E^C)$ is a diagram of spectra

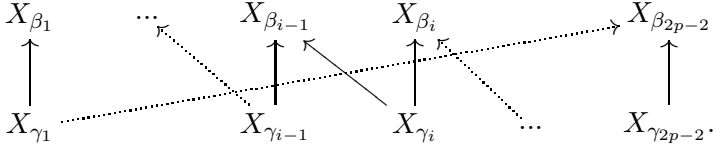
\[
\xymatrix{      X_{\beta_1} & ... & X_{\beta_{i-1}} & X_{\beta_i} & & X_{\beta_{2p-2}} \\
             X_{\gamma_1} \ar[u]\ar@{.>}[urrrrr] & & X_{\gamma_{i-1}}\ar[u]\ar@{.>}[ul] & X_{\gamma_i}\ar[u]\ar[ul] & ... \ar@{.>}[ul]
            & X_{\gamma_{2p-2}}. \ar[u]
}\]

{\bf N.B.} It should be pointed out that we work in the homotopy category of a diagram category of spectra and not
with diagrams taking values in the homotopy category of spectra.

In this particular case it it not too hard to characterise the fibrant and cofibrant objects
of $\E^C$:

\begin{itemize}
\item $X \in \E^C$ is fibrant iff each $X_{\beta_i}, X_{\gamma_i}$ is fibrant in $\E$
\item $X \in \E^C$ is cofibrant iff each $X_{\beta_i}, X_{\gamma_i}$ is cofibrant in $\E$
and for all $i \in \mathbb{Z} \big/ {\scriptstyle (2p-2)}$.
\[
X_{\gamma_{i+1}} \vee X_{\gamma_i} \longrightarrow X_{\beta_i}
\]
is a cofibration in $\E$.
\end{itemize}

\section{The functor $\Q$}

\subsection{Defining $\Q$}

We would now like to build twisted cochain complexes out of diagrams of spectra. Let $X$
be an object of $\Ho(\E^C)$. The given morphism
\[
p_i: X_{\gamma_{i}} \longrightarrow X_{\beta_i}
\]
as a part of the diagram $X$ induces a morphism in $\A$

\[
\pi_i:= E(1)_*(p_i)[i]: E(1)_*(X_{\gamma_i})[i] \longrightarrow E(1)_*(X_{\beta_i})[i].
\]

\begin{nnn}
$G^i(X):= E(1)_*(X_{\gamma_i})[i]$ and $B^i(X):=E(1)_*(X_{\beta_i})[i] $.
\end{nnn}

The objects $B^i(X)$ will play the role of the boundaries in the
cochain complex $C^*(X)$ to be built, and the $G^i(X)$'s will play
the role of the quotient of the cochains by the boundaries.

Now we would like to assign to each $k_i: X_{\gamma_{i+1}} \longrightarrow X_{\beta_i} \in \Ho(\E^{\underline{1}})$ (see section 1.5) an exact triangle
\[
X_{\gamma_{i+1}} \stackrel{k_i}{\longrightarrow} X_{\beta_i} \longrightarrow \cone(k_i)
\longrightarrow \Sigma X_{\gamma_{i+1}}.
\]
in a functorial (!) way. This is done by using Franke's cone functor
\[
\cone: \Ho(\E^{\underline{1}}) \longrightarrow \Ho(\E), \, \, \, (f:A \rightarrow B) \mapsto \Hocolim(\ast\leftarrow A \stackrel{f}{\rightarrow}B).
\]

\begin{nnn}
Define $C^i(X):= E(1)_*(\cone(k_i))[i] \in \A$.
\end{nnn}

Applying $E(1)_*$ to the above exact triangle we obtain a long exact sequence

\begin{equation}
...\rightarrow G^{i+1}(X)[-1]\rightarrow B^i(X) \rightarrow C^i(X) \rightarrow G^{i+1}(X)
\rightarrow B^i(X)[1] \rightarrow...
\end{equation}

Now let $\mathcal{L}$ be the full subcategory of $\Ho(\E^C)$ consisting of all objects $X$
such that

\begin{itemize}
\item $G^i(X)$ and $B^i(X)$ are not just objects of $\A$ but actually objects of the splitting $\B$ of $\A$ (see section 1.1).
\item $\pi_i: G^i(X) \longrightarrow B^i(X)$ is surjective for all $i$.
\end{itemize}

So if $X$ is an object of $\mathcal{L}$, what does this mean for the long exact
sequence (1)? If $X \in \mathcal{L}$, then by definition
\[
G^{i+1}(X)[-1] \in \B[-1] \quad\mbox{and}\quad B^i(X) \in \B.
\]
Therefore, by definition of $\B$, the maps $G^{i+1}(X)[-1] \longrightarrow B^i(X)$ and \linebreak
$G^{i+1}(X) \longrightarrow B^i(X)[1]$ in the long exact sequence (1) are zero.
Thus, (1) splits into short exact sequences
\begin{equation}
0 \longrightarrow B^i(X) \stackrel{\iota_i}{\longrightarrow} C^i(X) \stackrel{\rho_i}{\longrightarrow}
G^{i+1}(X) \longrightarrow 0.
\end{equation}

To make a cochain complex out of the objects $C^i(X)$, we need a differential \linebreak $d:C^i(X) \longrightarrow
C^{i+1}(X)$ which we define as the composition

\begin{equation}
C^i(X) \stackrel{\rho_i}{\longrightarrow} G^{i+1}(X) \stackrel{\pi_{i+1}}{\longrightarrow}
B^{i+1}(X) \stackrel{\iota_{i+1}}{\longrightarrow} C^{i+1}(X).
\end{equation}

Then $d^2$ is zero indeed since it factors over $\rho_{i+1} \circ \iota_{i+1}$ which is part
of the short exact sequence (2) and thus zero itself. The morphisms $\rho_i$ and $\pi_i$ are surjective
since $X \in \mathcal{L}$, so $\im(d)=B^*(X)$. Also, because of the shape of the underlying
poset we work with, $C^*(X)$ is $2p-2$-twistperiodic. So this construction
gives a functor

\[
\hspace{3cm}\Q: \mathcal{L} \longrightarrow \Cb, \quad X \longmapsto C^*(X).
\]

The next aim is to show that $\Q$ is an equivalence of categories which will be done
in the next two subsections.

\vsp
\subsection{$\Q$ is full and faithful}

We have to prove that for objects $X$ and $\xtild$ of $\mathcal{L}$, the map

\begin{equation}
M:=\Hom_{\Ho(\E^C)}(X,\xtild) \stackrel{q}{\longrightarrow} N
\end{equation}

with
\[
N:=\bigoplus\limits_{i}\Hom_{\B^{\underline{1}}}\big(((B^i(X)\rightarrow C^i(\xtild)),(B^i(\xtild) \rightarrow C^i(\xtild))\big)
\]

induced by $\Q$ is injective and its image consists of those morphisms that are morphisms
of cochain complexes. A morphism $f=(f_i)_i \in N$ is also a morphism
of cochain complexes iff it is compatible with the differentials, i.e. (remembering the
definition of $d$) makes the outer square in the following diagram commute:

\[
\xymatrix{ C^i(X) \ar[r]^{\rho_i}\ar[d]_{f^i} & G^{i+1}(X) \ar[r]^{\pi_{i+1}} \ar[d]_{\overline{f}^{i}}
             & B^{i+1}(X) \ar[r]^{\iota_{i+1}}\ar[d]^{f^{i+1}} & C^{i+1}(X) \ar[d]^{f^{i+1}} \\
            C^i(\xtild) \ar[r]^{\rho_i} & G^{i+1}(\xtild) \ar[r]^{\pi_{i+1}} & B^{i+1}(\xtild) \ar[r]^{\iota_{i+1}} & C^{i+1}(\xtild)
}
\]

Since $f \in N$ and $G^{i+1} \cong C^i \big/ B^i$, we know that the first and the third
small square commute. So, $f$ is a morphism of cochain complexes if and only if the middle
small square commutes, i.e. iff $f$ lies in the kernel of the map

\[
D: N \longrightarrow \bigoplus\limits_i \Hom_{\A}\big( G^i(X), B^i(\xtild) \big)
\]

where $D$ sends $f=(f_i)_i \in N$ to $f^{i+1}\circ \pi_{i+1} - \pi_{i+1}\circ\overline{f}^i$, with  $\overline{f}^i: G^{i+1}(X) \rightarrow G^{i+1}(\xtild)$ induced by $f^i$.

So, showing that $\Q$ is full and faithful is equivalent to showing that

\begin{equation}
0 \longrightarrow M \stackrel{q}{\longrightarrow} N \stackrel{D}{\longrightarrow}
\bigoplus\limits_i \Hom_{\A}\big( G^i(X), B^i(\xtild) \big)
\end{equation}

is exact. To show the exactness of (5), we would first like to get a description of $M$ and $N$ in
terms of some other exact sequences.

\vsp We start with $M$. A morphism of $\Hom_{\Ho(\E^C)}(X,\xtild)$
consists of the following data: the morphisms at each vertex plus
commutativity conditions coming from the shape of $C$. To be more
precise, the mapping space $\map_{\E^C}(X,\xtild)$ (see Section
\ref{quillen}) is the upper left corner of the following pullback
square of mapping spaces

\[
\xymatrix{ \map_{\E^C}(X,\xtild) \ar[r] \ar[d] & \prod\limits_i \map_{\E}(X_{\beta_i}, \xtild_{\beta_i}) \ar[d] \\
           \prod\limits_i \map_{\E}(X_{\gamma_i}, \xtild_{\gamma_i}) \ar[r] & \prod\limits_i \map_{\E}(X_{\gamma_{i+1}}, \xtild_{\beta_i}) \times \prod\limits_i \map_{\E}(X_{\gamma_i}, \xtild_{\beta_i})
}
\]

where the lower left and upper right corner contain the information about the maps at each
vertex and the lower right corner plus the maps into it give the commutativity conditions.
The right vertical map is the precomposition with the maps
\begin{equation}
X_{\gamma_{i+1}} \vee X_{\gamma_i} \longrightarrow X_{\beta_i}
\end{equation}
and the lower horizontal map is the composition with the maps
\[
\xtild_{\gamma_i} \longrightarrow \xtild_{\beta_i} \quad\quad \mbox{resp.} \quad\quad
\xtild_{\gamma_i} \longrightarrow \xtild_{\beta_{i-1}}.
\]
Without loss of generality one can assume $X$ to be cofibrant and
$\xtild$ to be fibrant (see section 1.3). Since (6) is then a
cofibration for each $i$ and $\E$ is a simplicial model category
(see e.g. \cite{GJ99}, section II.3), the right vertical map in the
pullback square is a fibration. Therfore, the pullback square is a
homotopy pullback square, and the left vertical map is a fibration
as well.

From a homotopy pullback square one gets a long exact homotopy sequence. Since $X$ is cofibrant
and $\xtild$ fibrant, we have as homotopy groups
\[
\pi_k \map_{\E}(X_{\gamma_i}, \xtild_{\gamma_i})\cong [X_{\gamma_i}, \xtild_{\gamma_i}]_{k}^{E(1)}
\]
(analogously for the other indices), and
\[
\pi_0 \map_{\E^C}(X,\xtild) = M \cong \Hom_{\Ho(\E^C)}(X,\xtild).
\]
Here, $[A,B]^{E(1)}_k$ denotes $\Hom_{\Ho(\E)}(\Sigma^k A,B)$.
Writing down the first five terms of the long exact homotopy sequence we obtain

\begin{equation}
\xymatrix{
\bigoplus\limits_i [X_{\gamma_i},\xtild_{\gamma_i}]_{1}^{E(1)}  \oplus
\bigoplus\limits_i [X_{\beta_i},\xtild_{\beta_i}]_{1}^{E(1)}  \ar[d]
\\
\bigoplus\limits_i [X_{\gamma_{i+1}},\xtild_{\beta_i}]_{1}^{E(1)}  \oplus
\bigoplus\limits_i [X_{\gamma_i},\xtild_{\beta_i}]_{1}^{E(1)} \ar[d]
\\
{\textstyle M} \ar[d]
\\
\bigoplus\limits_i [X_{\gamma_i},\xtild_{\gamma_i}]_0^{E(1)}  \oplus
\bigoplus\limits_i [X_{\beta_i},\xtild_{\beta_i}]_0^{E(1)} \ar[d]
\\
 \bigoplus\limits_i [X_{\gamma_{i+1}},\xtild_{\beta_i}]_0^{E(1)}  \oplus
\bigoplus\limits_i [X_{\gamma_i},\xtild_{\beta_i}]_0^{E(1)}
}
\end{equation}

Next, we would like to simplify the terms of this sequence with the help of the $E(1)$-
Adams spectral sequence

\begin{equation}
E^{s,t}_2 \cong \Ext^{s}_{\A}(E(1)_{*+t}(Y),E(1)_*(Z)) \Rightarrow [Y,Z]_{t-s}^{E(1)}
\end{equation}

for $Y, Z \in \E$. Since in our case $X, \xtild \in \mathcal{L}$, we have
\[
E(1)_*(X_{\beta_i}), E(1)_*(\xtild_{\beta_i}) \in \B[-i].
\]
It follows that
\[ \Ext_{\A}^{s}(E(1)_{*+t}(X_{\beta_i}), E(1)_*(\xtild_{\beta_i}))\] is actually
$\Ext_{\B}^{s}(E(1)_{*+t}(X_{\beta_i}), E(1)_*(\xtild_{\beta_i}))$,
and by definition of $\B$, this $\Ext$-term can only be nonzero if
$t$ is a multiple of $2p-2$. This is because for an object of $\B$,
all objects in a injective resolutions are in $\B$ themselves again.
Bousfield also proved that $\Ext_{\A}^s(-,-)=0$ for $s \ge 3$.
Consequently, the spectral sequence collapses, as seen in the
following picture of the $E_2$-term for $p=3$:

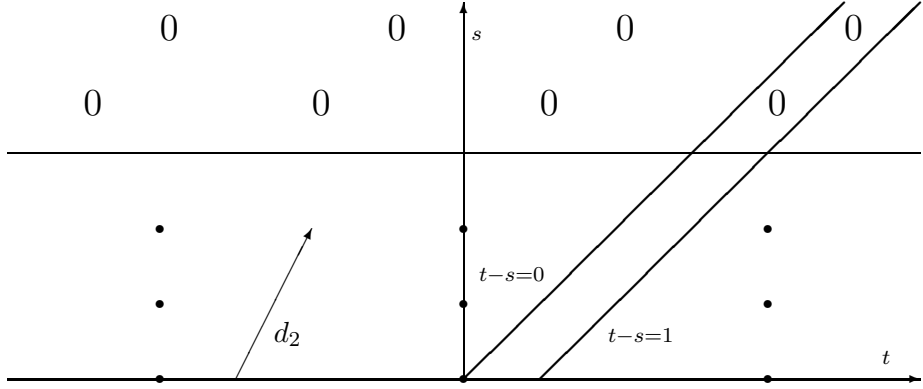
\setlength{\unitlength}{1cm}
\begin{picture}(25,7)

\put(0,0){\vector(1,0){12}}
\put(6,0){\vector(0,1){5}}
\put(6,0){\circle*{0.1}}
\put(6,1){\circle*{0.1}}
\put(6,2){\circle*{0.1}}
\put(2,0){\circle*{0.1}}
\put(2,1){\circle*{0.1}}
\put(2,2){\circle*{0.1}}
\put(10,0){\circle*{0.1}}
\put(10,1){\circle*{0.1}}
\put(10,2){\circle*{0.1}}

\put(1,3.5){\Large 0}
\put(4,3.5){\Large 0}
\put(7,3.5){\Large 0}
\put(10,3.5){\Large 0}

\put(2,4.5){\Large 0}
\put(5,4.5){\Large 0}
\put(8,4.5){\Large 0}
\put(11,4.5){\Large 0}

\put(3.5,0.5){$d_2$}
\put(3,0){\vector(1,2){1}}

\put(6.2,1.3){$\scriptstyle t-s=0$}
\put(7.9,0.5){$\scriptstyle t-s=1$}

\put(6.1,4.5){$\scriptstyle s$}
\put(11.5,0.2){$\scriptstyle t$}

\thinlines
\put(0,3){\line(1,0){12}}

\thicklines
\put(6,0){\line(1,1){5}}
\put(7,0){\line(1,1){5}}

\end{picture}

\vsp The $E_2$-term can only be nonzero at the location of the dots.
In particular, as this picture indicates, for all odd primes,
$E^{s,t}_2$ is zero if $t=s, s \neq 0$ and $t-s=1$. Therefore,
\[
[X_{\beta_i},\xtild_{\beta_i}]_{1}^{E(1)}=0=[X_{\gamma_i},\xtild_{\gamma_i}]_{1}^{E(1)}=[X_{\gamma_i},\xtild_{\beta_i}]_{1}^{E(1)}
\]

and

\begin{eqnarray*}
[X_{\beta_i},\xtild_{\beta_i}]_0^{E(1)} \cong \Hom_{\B}(E(1)_*(X_{\beta_i}),E(1)_*(\xtild_{\beta_i}))  \\
{[}X_{\gamma_i},\xtild_{\gamma_i}{]}_0^{E(1)} \cong \Hom_{\B}(E(1)_*(X_{\gamma_i}),E(1)_*(\xtild_{\gamma_i})) \\
{[}X_{\gamma_i},\xtild_{\beta_i}{]}_0^{E(1)} \cong \Hom_{\B}(E(1)_*(X_{\gamma_i}),E(1)_*(\xtild_{\beta_i}))
\end{eqnarray*}

Similarly, $\Ext_{\A}^{s}(E(1)_{*+t}(X_{\gamma_{i+1}}), E(1)_*(\xtild_{\beta_i}))$ can only
be non-zero if $s \le 2$ and \linebreak $t \equiv 1 (2p-2)$, in particular it is zero for $t-s=1, s \neq 0$
and $s=t, s \neq 1$.
So this spectral sequence also collapses, and it follows that

\[
[X_{\gamma_{i+1}},\xtild_{\beta_i}]_{1}^{E(1)} \cong \Hom_{\B}(E(1)_{*+1}(X_{\gamma_{i+1}}),E(1)_*(\xtild_{\beta_i}))
\]
and
\[
[X_{\gamma_{i+1}},\xtild_{\beta_i}]_0^{E(1)} \cong \Ext^1_{\B}(E(1)_*(X_{\gamma_{i+1}}),E(1)_*(\xtild_{\beta_i})).
\]

\pagebreak
Putting this into the sequence (7), we obtain the exact sequence

\begin{equation}
\xymatrix{
0 \ar[d] \\
\bigoplus\limits_i \Hom_{\B}(G^{i+1}(X),B^i(\xtild)) \ar[d] \\
M \ar[d] \\
\bigoplus\limits_i \Hom_{\B}(G^{i}(X),G^i(\xtild))
\oplus
\bigoplus\limits_i \Hom_{\B}(B^{i}(X),B^i(\xtild)) \ar[d] \\
\bigoplus\limits_i \Ext^1_{\B}(G^{i+1}(X),B^i(\xtild))
\oplus
\bigoplus\limits_i \Hom_{\B}(G^{i}(X),B^i(\xtild)).
}
\end{equation}

Now we would like to find a similar description of
\[
N=\bigoplus\limits_{i}\Hom_{\B^{\underline{1}}}\big(((B^i(X)\rightarrow C^i(X)),(B^i(\xtild) \rightarrow C^i(\xtild))\big).
\]
As mentioned before, morphisms in $N$ can be viewed as morphisms of the short exact sequences
\[
\xymatrix{ 0 \ar[r] & B^i(X) \ar[r]\ar[d]^{f_i} & C^i(X) \ar[r]\ar[d]^{f_i} & G^{i+1}(X) \ar[r]\ar[d]^{\overline{f}_i} & 0 \\
 0 \ar[r] & B^i(\xtild) \ar[r] & C^i(\xtild) \ar[r] & G^{i+1}(\xtild) \ar[r] & 0.
}
\]

Thus, we get a canonical map

\begin{equation}
N \longrightarrow N' := \begin{array}{c}
\bigoplus\limits_i \Hom_{\B}(B^i(X), B^i(\xtild)) \\
\oplus \\
\bigoplus\limits_i \Hom_{\B}(G^i(X), G^i(\xtild))
\end{array}
\end{equation}

by sending $f \in N$ to $ (f_i, \overline{f}_i)_i$. The kernel of this map consists
of morphisms of the same exact sequences of the form

\[
\xymatrix{ 0 \ar[r] & B^i(X) \ar[r]\ar[d]^{0} & C^i(X) \ar[r]\ar[d]^{\Phi} & G^{i+1}(X) \ar[r]\ar[d]^{0} & 0 \\
 0 \ar[r] & B^i(\xtild) \ar[r] & C^i(\xtild) \ar[r] & G^{i+1}(\xtild) \ar[r] & 0.
}
\]
Every $\Phi$ of the form
\[
C^i(X) \longrightarrow G^{i+1}(X) \stackrel{\phi}{\longrightarrow} B^i(\xtild) \longrightarrow C^i(\xtild)
\]
lies in the kernel of (10). From applying the snake lemma to the above diagram it also follows that every $\Phi$ in the kernel looks exactly like this. Therefore, the
kernel of (10) is isomorphic to $ \bigoplus\limits_i \Hom_{\B}(G^{i+1}(X),B^i(\xtild))$.
Consequently,

\begin{equation}
0 \longrightarrow \bigoplus\limits_i \Hom_{\B}(G^{i+1}(X),B^i(\xtild)) \longrightarrow N \longrightarrow N'
\end{equation}

is exact.

The next question is: when is an element of $N'$ hit by an element of $N$? In other words,
given $f_B: B^i(X) \rightarrow B^i(\xtild)$ and $f_G: G^{i+1}(X) \rightarrow G^{i+1}(\xtild)$,
when is there a map $f_C: C^i(X) \rightarrow C^i(\xtild)$ making the following diagram commute?

\[
\xymatrix{ 0 \ar[r] & B^i(X) \ar[r]\ar[d]^{f_B} & C^i(X) \ar[r]\ar[d]^{f_C} & G^{i+1}(X) \ar[r]\ar[d]^{f_G} & 0 \\
 0 \ar[r] & B^i(\xtild) \ar[r] & C^i(\xtild) \ar[r] & G^{i+1}(\xtild) \ar[r] & 0
}
\]

The upper sequence corresponds to an element $S \in \Ext^1_{\B}(G^{i+1}(X),B^i(X))$, the lower
one to an element $\tilde{S} \in \Ext^1_{\B}(G^{i+1}(\xtild),B^i(\xtild))$. The maps $f_B$ and $f_G$
give rise to maps

\begin{eqnarray*}
(f_B)_* : \Ext^1_{\B}(G^{i+1}(X),B^i(X)) \longrightarrow \Ext^1_{\B}(G^{i+1}(X),B^i(\xtild)) \\
(f_G)^* : \Ext^1_{\B}(G^{i+1}(\xtild),B^i(\xtild)) \longrightarrow \Ext^1_{\B}(G^{i+1}(X),B^i(\xtild)).
\end{eqnarray*}

So for given $f_B$ and $f_G$ there is a morphism $f_C$ making the above diagram commute if
and only if $(f_B)_*(S) = (f_G)^*(\tilde{S})$. \newpage It follows that

\begin{equation}
\xymatrix{
0 \ar[d] \\ \bigoplus\limits_i \Hom_{\B}(G^{i+1}(X),B^i(\xtild)) \ar[d] \\ N \ar[d] \\
N'=\bigoplus\limits_i \Hom_{\B}(B^i(X), B^i(\xtild))
\oplus
\bigoplus\limits_i \Hom_{\B}(G^i(X), G^i(\xtild))
 \ar[d]
\\ \bigoplus\limits_i \Ext^1_{\B}(G^{i+1}(X),B^i(\xtild))
}
\end{equation}
 is exact where the last map sends a pair of tupels $(f_B,f_G)$ to $(f_B)_*(S) - (f_G)^*(\tilde{S})$.
Putting this sequence together with the sequence (9), we obtain

\[
\xymatrix{
0 \ar[d] & 0 \ar[d] \\
\bigoplus\limits_i \Hom_{\B}(G^{i+1}(X),B^i(\xtild)) \ar[d]_a
\ar@{=}[r] &
\bigoplus\limits_i \Hom_{\B}(G^{i+1}(X),B^i(\xtild)) \ar[d]^b \\
M\ar[d]\ar[r]^q & N \ar[d] \\
{
\begin{array}{c}
\bigoplus\limits_i \Hom_{\B}(B^i(X), B^i(\xtild)) \\
\oplus \\
\bigoplus\limits_i \Hom_{\B}(G^i(X), G^i(\xtild))
\end{array}
}
 \ar@{=}[r]\ar[d] &
{
\begin{array}{c}
\bigoplus\limits_i \Hom_{\B}(B^i(X), B^i(\xtild)) \\
\oplus \\
\bigoplus\limits_i \Hom_{\B}(G^i(X), G^i(\xtild))
\end{array}
}
 \ar[d] \\
{\begin{array}{c}
\bigoplus\limits_i \Ext^1_{\B}(G^{i+1}(X),B^i(\xtild)) \\
 \oplus \\ \bigoplus\limits_i \Hom_{\B}(G^{i}(X),B^i(\xtild))
\end{array}} \ar[r]^{pr} &
\bigoplus\limits_i \Ext^1_{\B}(G^{i+1}(X),B^i(\xtild))
}
\]

where the second horizontal arrow is the morphism
induced by the functor $\Q$ and the last one is the projection onto the first summand. One
has to check that all the squares actually commute, which they do.

Then, a small diagram chase shows that $q$ is injective. Also, by
construction of $q$, in

\begin{equation}
0 \longrightarrow M \stackrel{q}{\longrightarrow} N \stackrel{D}{\longrightarrow}
\bigoplus\limits_i \Hom_{\A}\big( G^i(X), B^i(\xtild) \big),
\end{equation}

the image of $q$ lies in the kernel of $D$. With a slightly bigger diagram chase it follows
that the image of $q$ is the entire kernel of $D$.

This completes the proof that $\Q$ is full and faithful.

\subsection{$\Q$ is essentially surjective}

To complete the proof of the claim that
\[
\Q: \mathcal{L} \longrightarrow \Cb
\]
is an equivalence of categories, it is left to
show that $\Q$ is essentially surjective, i.e. each $C^* \in \Cb$ is isomorphic to an
object the image of $\Q$. So let $C^*$ be an object of $\Cb \cong \Ca$, and let
$B^*(C)$ denote the boundaries of $C^*$ and $G^*(C)$ the quotient of $C^*$ by its boundaries.
We will prove our claim by induction on the injective dimension of the $B^i(C)$'s and
$G^i(C)$'s. That means, we will perform an induction on $k$ where

\[
\ma_i(\injdim B^i(C), \injdim G^i(C)) \le k \le 2.
\]

Let $I \in \A$ be an injective object, and consider the following cochain complexes:

\[
V(I)^* \quad\mbox{with}\quad V(I)^n:= T^{n(p-1)}(I), \,\, d=0
\]
\[
C(I)^* \quad\mbox{with}\quad C(I)^n:= T^{n(p-1)}(I) \oplus T^{(n-1)(p-1)}(I),
\quad d=  \left( \begin{array}{cc} 0 & 0 \\
id & 0
\end{array} \right)
\]

with the structure isomorphisms $\alpha_{V(I)}$ and $\alpha_{C(I)}$ (see 1.2) being the identity.
Both complexes are injective in $\Ca \cong \Cb$. The complex $V(I)^*$ belongs to the essential image
of $\Q$: Without loss of generality, let $I$ be an object of $\B$.
First, this complex can be realized by spectra $X_i \in \E$ such that
\[
E(1)_*(X_i)[i] \cong T^{i(p-1)}(I),
\]
see e.g. \cite{Bou85}, Prop. 8.2. Now look at the following diagram
$X \in \Ho(\E^C)$:

\[
\xymatrix{     \ast  & & \ast & \ast & & \ast \\
             \Sigma^{-1} X_{1} \ar[u]\ar@{.>}[urrrrr] & & \Sigma^{-1} X_{i-1}\ar[u]\ar@{.>}[ul] &
          \Sigma^{-1} X_{i} \ar[u]\ar[ul] &\ar@{.>}[ul] & \Sigma^{-1} X_{2p-2} \ar[u]
}
\]

with $X_i$ as above. Clearly, $X \in \mathcal{L}$. Applying
$Q$ to this diagram $X$ one sees that

\[
C^i(X) = E(1)_*(\cone(\Sigma^{-1} X_i \rightarrow \ast))[i] = E(1)_*(X_i)[i] \cong
T^{i(p-1)}(I) = V(I)^i
\]
together with the correct zero differential. This means that $V(I)^*$ is in the essential
image of $\Q$, and similarly, also $C(I)^*$.

Now, to begin the induction, let $C^*$ be a complex such that

\[
\ma_i(\injdim B^i(C), \injdim G^i(C)) = 0.
\]

It follows that $H^0(C)$ and $B^0(C)$ are injective objects of $\A$, and one checks that
\[
C^* \cong V(H^0(C))^* \oplus C(B^0(C))^*.
\]
Consequently, $C^*$ lies in the essential image of $\Q$ which starts the induction.

Next, let our claim be true for $k-1$ and let $C^*$ be an arbitrary complex with

\[
\ma_i(\injdim B^i(C), \injdim G^i(C)) \le k.
\]

$\Ca$ has enough injectives (\cite{Fra96}, Prop. 1.3.3), i.e. there is an embedding \linebreak
$C^* \stackrel{i}{\longrightarrow} K^*$ such that $K^*$ is strictly injective and $i$ is a
weak equivalence. Consequently,

\[
\ma_i(\injdim B^i(K), \injdim G^i(K)) =0.
\]

We have already proved that $K^*$ is in the essential image of $\Q$. Looking at

\begin{equation}
0 \longrightarrow C^* \stackrel{i}{\longrightarrow} K^* \longrightarrow L^*:= \coker(i)
\longrightarrow 0 ,
\end{equation}

we now prove that
\[
\ma_i(\injdim B^i(L), \injdim G^i(L)) \le k-1.
\]

For example, if $\injdim B^i(C) \le k$, then

\[
0 \longrightarrow B^i(C) \longrightarrow B^i(K) \longrightarrow B^i(L) \longrightarrow 0
\]
 is exact and $B^i(K)$ is injective. If

\[
B^i(L) \rightarrow J_0 \rightarrow J_1 \rightarrow ... \rightarrow J_m \rightarrow 0
\]

is an injective resolution of $B^i(L)$, then

\[
B^i(C) \rightarrow B^i(K) \rightarrow J_0 \rightarrow J_1 \rightarrow ... \rightarrow J_m \rightarrow 0
\]

is an injective resolution of $B^i(C)$. Since there is a resolution of $B^i(C)$ of length
$ \le k$, it follows that there is also a resolution for $B^i(L)$ of length $\le k-1$.

This shows that

\[
\ma_i(\injdim B^i(L), \injdim G^i(L)) \le k-1,
\]
and by our induction, $L^*$ lies in the essential image.

The fact that $C^*$ now lies in the essential image of $\Q$ as well is a consequence of the
following:

Let $Y_1 \rightarrow Y_2 \rightarrow Y_3 \rightarrow Y_1[1]$ be an exact triangle in
$\Ho(\E^C)$, and $Y_2, Y_3 \in \mathcal{L}$, \linebreak $\Q(Y_2) \rightarrow \Q(Y_3)$ an epimorphism
of cochain complexes and $H^*(Q(Y_2)) \rightarrow H^*(Q(Y_3))$ be surjective. Then
\[
0 \longrightarrow \Q(Y_1) \longrightarrow \Q(Y_2) \longrightarrow \Q(Y_3) \longrightarrow 0
\]
is exact and $Y_1$ is an object of $\mathcal{L}$.

(To prove this, one frequently uses the five lemma and has to remember that $\B$ is a Serre
class in $\A$.)

Back to our short exact sequence (14).
We have proved that there are objects \linebreak $X_2, X_3 \in \Ho(\E^C)$ such that $\Q(X_2) \cong K^*$ and
$\Q(X_3) \cong L^*$. Since we also know that $\Q$ is full, we see that the map
$\Q(X_2) \rightarrow \Q(X_3)$ is induced by a morphism $X_2 \rightarrow X_3 \in \Ho(\E^C)$
which can be completed to an exact triangle $X_1 \rightarrow X_2 \rightarrow X_3 \rightarrow
X_1[1]$. $\Q(X_2) \rightarrow \Q(X_3)$ is an epimorphism that also induces an epimorphism
in cohomology and thus satisfies the condition above. It follows that $X_1 \in
\mathcal{L}$, that
\[
0 \longrightarrow \Q(Y_1) \longrightarrow \Q(Y_2) \longrightarrow \Q(Y_3) \longrightarrow 0
\]
is exact and that therefore $C^* \cong \Q(X_1)$. This completes the proof that $\Q$ is essentially
surjective and consequently is an equivalence of categories.

\section{The reconstruction functor $\R$}

\subsection{Defining $\R$}

In the last section we showed that

\[
\hspace{3cm}\Q: \mathcal{L} \longrightarrow \Cb
\]

is an equivalence of categories. To prove the main theorem, we would like to build
an equivalence of categories

\[
\R: \Db = \Ho(\Cb) \longrightarrow \Ho(\E)
\]

with the help of $\Q$. Define

\[
\R' := \mbox{Hocolim} \circ \Q^{-1} : \Cb \longrightarrow \Ho(\E^C) \longrightarrow \Ho(\E).
\]

We would like to show that $\R'$ factors over the derived category of $\Cb$. This will gives
us the desired reconstruction functor $\R$ of which we would like to show that it is an
equivalence of categories.

However, we first look at some properties of
\[
E(1)_{*} \circ \R' : \Cb \longrightarrow \A.
\]

\begin{lem}\label{hocolim}
\[
E(1)_{*} (\mbox{Hocolim}_C X) \cong \bigoplus\limits_i H^i(\Q(X))[-i]
\]
\end{lem}

\begin{proof}

By definition,
\[
\mbox{Hocolim}_C X = \mbox{colim}_C X^{cof}
\]

where $X^{cof}$ denotes the cofibrant replacement of $X \in \E^C$. Now let us look at the
colimit of a diagram

\[
\xymatrix{      X_{\beta_1} & ... & X_{\beta_{i-1}} & X_{\beta_i} &   & X_{\beta_{2p-2}} \\
             X_{\gamma_1} \ar[u]\ar@{.>}[urrrrr] & & X_{\gamma_{i-1}}\ar[u]\ar@{.>}[ul] & X_{\gamma_i}\ar[u]\ar[ul] & ... \ar@{.>}[ul] & X_{\gamma_{2p-2}}. \ar[u]
}
\]

We have morphisms
\[
X_{\gamma_i} \vee X_{\gamma_{i+1}} \longrightarrow X_{\beta_i}
\]
for each $i$. Taking the wedge of those morphisms for even $i$, one obtains a morphism
\[
\bigvee\limits_{i=1}^{2p-2} X_{\gamma_i} \longrightarrow \bigvee\limits_{i \,\, \even} X_{\beta_i},
\]
and simultaneously, for odd $i$,
\[
\bigvee\limits_{i=1}^{2p-2} X_{\gamma_i} \longrightarrow \bigvee\limits_{i \,\, \odd} X_{\beta_i}.
\]
The colimit of the diagram $X$ is the same as the colimit of the following diagram:
\[
\bigvee\limits_{i \,\, \odd} X_{\beta_i} \longleftarrow \bigvee\limits_{i=1}^{2p-2} X_{\gamma_i} \longrightarrow \bigvee\limits_{i \,\, \even} X_{\beta_i},
\]
i.e. the colimit of $X$ is the pushout of the upper left corner in

\[
\xymatrix{ \bigvee\limits_{i=1}^{2p-2} X_{\gamma_i} \ar[d]\ar[r] & \bigvee\limits_{i \,\, \even} X_{\beta_i} \ar[d] \\
\bigvee\limits_{i \,\, \odd} X_{\beta_i} \ar[r] & \mbox{colim}_C X.
}
\]
Without loss of generality, let $X$ be cofibrant, so that the colimit of $X$ models the homotopy colimit.
Then the left vertical and upper horizontal maps in the square are cofibrations, and the pushout
diagram is also a homotopy pushout diagram. Therefore,
\[
\bigvee\limits_{i=1}^{2p-2} X_{\gamma_i} \rightarrow \bigvee\limits_{i \,\, \odd} X_{\beta_i} \vee \bigvee\limits_{i \,\, \even} X_{\beta_i} \cong \bigvee\limits_{i=0}^{2p-2} X_{\beta_i}
\rightarrow \mbox{Hocolim}_C X \rightarrow \Sigma\Big( \bigvee\limits_{i=1}^{2p-2} X_{\gamma_i} \Big)
\]
is an exact triangle in $\Ho(\E)$. Applying $E(1)$-homology, one obtains a long exact sequence
\begin{eqnarray}
...\bigoplus\limits_i E(1)_n(X_{\gamma_i}) \rightarrow \bigoplus\limits_i E(1)_n(X_{\beta_i})
\rightarrow E(1)_n(\mbox{Hocolim}_C X) \rightarrow \nonumber\\
\rightarrow \bigoplus\limits_i E(1)_{n-1}(X_{\gamma_i}) \rightarrow \bigoplus\limits_i E(1)_{n-1}(X_{\beta_i}) ...
\end{eqnarray}

The map
\[
\oplus \pi_i[-i+1]: \quad \bigoplus\limits_i E(1)_{n-1}(X_{\gamma_i}) \rightarrow \bigoplus\limits_i E(1)_{n-1}(X_{\beta_i})
\]
is surjective for all $n$ since $X \in \mathcal{L}$, so
\[
\bigoplus\limits_i E(1)_n(X_{\gamma_i}) \longrightarrow E(1)_n(\mbox{Hocolim}_C X)
\]
is the zero map. So we get a short exact sequence in $\A$
\[
0 \rightarrow E(1)_*(\mbox{Hocolim}_C X) \longrightarrow \bigoplus\limits_i E(1)_{*-1}(X_{\gamma_i})
\xrightarrow{\oplus\pi_i[-i+1]}  E(1)_{*-1}(X_{\beta_i}) \rightarrow 0.
\]

Therefore,
\[
E(1)_*(\mbox{Hocolim}_C X) \cong \bigoplus\limits_i \ke(\pi_i) [-i+1].
\]

Now we prove that $\ke(\pi_i)$ is isomorphic to $H^{i-1}(\Q(X))$.
Let us remember how the differential $d$
of $C^*(X) = \Q(X)$ had been defined (see section 2.1). Here is $d^2$:

\begin{eqnarray*}
C^i(X) \stackrel{\rho_i}{\rightarrow} G^{i+1}(X) \stackrel{\pi_{i+1}}{\rightarrow}
B^{i+1}(X) \stackrel{\iota_{i+1}}{\rightarrow} C^{i+1}(X) \stackrel{\rho_{i+1}}{\rightarrow} G^{i+2}(X) \stackrel{\pi_{i+2}}{\rightarrow} B^{i+2}(X) \stackrel{\iota_{i+2}}{\rightarrow} C^{i+2}(X)
\end{eqnarray*}

We have $\im(\iota_{i+1}) = \ke(\rho_{i+1})$ since they are part of the short exact sequence (2).
Since $\rho_i$ and $\pi_{i+1}$ are surjective, $\im(d) = \im(\iota_{i+1})$. We also have
$\ke(d)= \ke(\pi_{i+2}\circ\rho_{i+1})$. By basic algebra,

\[
\ke(\pi_{i+2}) \cong \frac{\ke(\pi_{i+2}\circ\rho_{i+1})}{\ke(\rho_{i+1})} \cong
\frac{\ke(d)}{\im(\iota_{i+1})} \cong \frac{\ke(d)}{\im(d)} \cong H^{i+1}(\Q(X)).
\]

It follows that

\[
E(1)_{*}(\mbox{Hocolim}_C X) \cong \bigoplus\limits_i H^i(\Q(X))[-i].
\]
\end{proof}

Because of the lemma we now see that the functor $E(1)_* \circ \R'$ sends weak equivalences
(i.e. quasi-isomorphisms) in $\Cb$ to isomorphisms in $\A$ and thus factors over $\Db = \Ho(\Cb)$.
In other words, for $C^*, D^*$ quasi-isomorphic cochain complexes we get
\begin{eqnarray*}
E(1)_{*}(\R'(C^*)) \cong \bigoplus\limits_i H^i(C^*)[-i] \cong \bigoplus\limits_i H^i(D^*)[-i]
\cong E(1)_{*}(\R'(D^*)).
\end{eqnarray*}

However, two objects of $\Ho(\E)$ are isomorphic if and only if there is a morphism of spectra
inducing an isomorphism in $E(1)$-homology,
so $\R'(C^*) \cong \R'(D^*)$ for quasi-isomorphic $C^*$ and $D^*$, and consequently $\R'$ itself factors over the derived category
$\Db$. So we have obtained a functor

\[
\R: \Db \longrightarrow \Ho(\E).
\]

\subsection{The main theorem}

\begin{ttt}
$\R$ is an equivalence of categories.
\end{ttt}

\begin{proof}
First again, we prove that $\R$ is full and faithful, i.e. for

\[C_1^*, C_2^* \in \Db \cong \Da,
\]
the map
\[
r: \Hom_{\Da}(C_1^*, C_2^*) \longrightarrow [\R(C_1^*), \R(C_2^*)]^{E(1)}
\]
induced by $\R$ is an isomorphism.

To show this, we once more make use of the Adams spectral sequence (\cite{Fra96} 2.1.1)

\begin{equation}
E_2^{s,t} = \Ext^s_{\A}\Big(\bigoplus\limits_i H^i(C_1^*)[-i-t],\bigoplus\limits_i H^i(C_2^*)[-i]\Big) \Rightarrow \Hom_{\Da}(C_1^*[t-s],C_2^*)
\end{equation}

where $C_1^*, C_2^* \in \Da$. This spectral
sequence arises as follows: We begin with an injective resolution of $\bigoplus\limits_i H^i(C_2^*)[-i]$:

\begin{equation}
\label{injres}
\xymatrix{ \bigoplus H^i(C_2^*)[-i] \ar@{^{(}->}[r] & I^0 \ar[r]^{d_1} \ar@{>>}[d]  & I^1 \ar[r]^{d_2} \ar@{>>}[d] & I^2 \ar[r] & 0 \\
 & im(d_1) \ar@{^{(}->}[ur] & im(d_2) \ar@{^{(}->}[ur] & &
}
\end{equation}

(This resolution stops at $I^2$ since the injective dimension of an object in $\A$ is at most 2.)

This resolution gives rise to an Adams resolution

\begin{equation}
\label{adamsres1}
\xymatrix{
C_2^* = C_2^{(0)} \ar[d] & C_2^{(1)} \ar[l]\ar[d] & C_2^{(2)} \ar[l]\ar[d] & 0 \ar[l] \\
E_{I^0} \ar[ur]_{+} & E_{I^1} \ar[ur]_{+} & E_{I^2} \ar[ur]_{+} &
}
\end{equation}

The Adams resolution is characterised by the following: First, by applying
\[
\bigoplus\limits_i H^i(-)[-i]
\]

to the diagram

\begin{equation}
\label{adamsres2}
\xymatrix{ C_2^* = C_2^{(0)} \ar[r] & E_{I^0} \ar[r] \ar[d]  & E_{I^1} \ar[r] \ar[d] & E_{I^2} \ar[r] & 0 \\
 & C_2^{(1)} \ar[ur] & C_2^{(2)} \ar[ur] & &
}
\end{equation}

one obtains exactly the diagram (\ref{injres}).
Besides, each triangle in (\ref{adamsres1})
is an exact triangle in $\Da$ (the diagonal maps are maps raising the degree by one), and $E_I$ denotes the Eilenberg-MacLane object
for $I \in \A$, i.e.
\[
\Hom_{\A}(\bigoplus\limits_i H^i(C^*)[-i],I) \cong \Hom_{\Da}(C^*,E_I) \,\, \mbox{for all} \,\, C^* \in \Da,
\]
and for $C^* = E_I$, the image of the identity in
\[
\Hom_{\A}(\bigoplus\limits_i H^i (E_I)[-i],I)
\]
is an isomorphism. (Note that by Lemma 2.1.1 of \cite{Fra96},
$C_2^{(2)}$ is an Eilenberg-MacLane object for $I^2$ indeed!)
Applying $\Hom_{\Da}(C_1^*,-)$ to the resolution (\ref{adamsres1})
gives an exact couple, and with it the desired spectral sequence.

We now apply the reconstruction functor $\R$ to (\ref{adamsres1}) and claim that the result
\begin{equation}
\label{Radamsres}
\xymatrix{ \R(C_2) = \R(C_2^{(0)}) \ar[d] & \R(C_2^{(1)}) \ar[l]\ar[d] & \R(C_2^{(2)}) \ar[l] \ar[d]_{\cong} & 0 \ar[l] \\
\R(E_{I^0}) \ar[ur]_{+} & \R(E_{I^1}) \ar[ur]_{+} & \R(E_{I^2}) \ar[ur]_{+} &
}
\end{equation}

is an Adams resolution for $\R(C^*_2)$ with respect to $E(1)$-homology.

\vsp
We have to prove:
\begin{itemize}
\item applying $E(1)_*$ to (\ref{Radamsres}) gives an injective resolution of $E(1)_*(\R(C_2^*))$
\item each triangle in (\ref{Radamsres}) is exact
\item $\R(E_I)$ is again an Eilenberg-MacLane object in $\Ho(\E)$
\end{itemize}

\vsp
The first point is clear after the Lemma \ref{hocolim}, which says that
\[
E(1)_*(\R(C^*)) \cong \bigoplus\limits_i H^i(C^*)[-i].
\]
To prove the second point we make use of the following
fact without giving the details of its proof:

Let $C_0^* \rightarrow C_1^* \rightarrow C_2^* \rightarrow C_0^*[1]$ be an exact triangle
in $\Da$ with $H^*(C_0^*) \rightarrow H^*(C_1^*)$ a monomorphism. Then
\[
\R(C_0^*) \rightarrow \R(C_1^*) \rightarrow \R(C_2^*) \rightarrow \R(C_0^*[1])
\]
is an exact triangle in $\Ho(\E)$.

\vsp
Using the Lemma \ref{hocolim} again,
we see that the vertical arrows in (\ref{adamsres1}) give monomorphisms in cohomology.
So, applying the above fact, we have that the triangles in (\ref{adamsres1}) are exact indeed.

\vsp
To show that $\R(E_I)$ is again an Eilenberg-MacLane object in $\Ho(\E)$ for injective $I \in \A$,
we have to show that
\[
\Hom_{\A}(E(1)_*(X),I) \cong [X,\R(E_I)]^{E(1)} \quad\mbox{for all} \quad X \in \Ho(\E).
\]
We know that
\[
E(1)_*(\R(E_I)) \cong \bigoplus\limits_i H^i (E_I) [-i] \cong I,
\]
so $\R(E_I)$ has injective $E(1)$-homology. Now we look at the classical Adams spectral sequence

\begin{eqnarray*}
E_2^{s,t} = \Ext_{\A}^s(E(1)_*(X),E(1)_*(\R(E_I))[t]) = \Ext_{\A}^s(E(1)_*(X),I[t]) \\
\Rightarrow [X,\R(E_I)]^{E(1)}_{t-s}
\end{eqnarray*}

for $X \in \Ho(\E)$.
Since $I$ is injective in $A$, the $\Ext$-term vanishes unless $s=0$, so the spectral sequence
collapses, and
\[
\Ext_{\A}^0(E(1)_*(X),I[t]) = \Hom_{\A}(E(1)_*(X),I[t]) \cong [X,\R(E_I)]^{E(1)}_{t}
\]
as desired.

\vsp Applying $[\R(C_1^*),-]^{E(1)}$ to (\ref{Radamsres}) gives an
exact couple leading to the Adams spectral sequence

\[
E_2^{s,t}=\Ext_{\A}^s(E(1)_*(\R(C_1^*)),E(1)_*(\R(C_2^*))[t])
 \Rightarrow [\R(C_1^*),\R(C_2^*))]^{E(1)}_{t-s}.
\]

So $\R$ induces a morphism of exact couples that is also an isomorphism on the $E_1$-terms

\[
r: \Hom_{\Da}^t(C_1^*,E_{I^s}) \longrightarrow [\R(C_1^*), \R(E_{I^s})]^{E(1)}_t:
\]
by definition of an Eilenberg-MacLane object, the left side is
isomorphic to $$\Hom^t_{\A}( \bigoplus\limits_i H^i (C_1^*) [-i],
I^s).$$  Since $\R(E_{I^s})$ is an Eilenberg-MacLane object with
respect to $E(1)_*$, the right side is isomorphic to
$$\Hom_{\A}^t(E(1)_*(\R(C^*_1),I^s).$$ So because of Lemma
\ref{hocolim}, both sides are isomorphic. It follows that $r$ is an
isomorphism on the targets of the spectral sequences, and thus, $\R$
is full and faithful.

\vsp
Now it is left to show that $\R$ is essentially surjective. Let $Y$ be an object of $\Ho(\E)$
and let
\begin{equation}
\label{Yadamsres}
\xymatrix{ Y = Y^{(0)} \ar[d] & Y^{(1)} \ar[l]\ar[d] & Y^{(2)} \ar[l]\ar[d]_{\cong} & 0 \ar[l] \\
\mathcal{E}_{I^0} \ar[ur]_{+} & \mathcal{E}_{I^1} \ar[ur]_{+} & \mathcal{E}_{I^2} \ar[ur]_{+} &
}
\end{equation}
be an Adams resolution for $Y$. First, we show that all Eilenberg MacLane objects \linebreak $\mathcal{E}_I \in \Ho(\E)$ lie in
the essential image of $\R$: Let $E_I$ be the Eilenberg-MacLane object for $I$ in $\Da$. We
already showed that $\R(E_I)$ is an Eilenberg-MacLane object for $I$ in $\Ho(\E)$, and thus,
$\mathcal{E}_I \cong \R(E_I)$, so $\mathcal{E}_I$ lies in the essential image of $\R$.

Next, we would like to show that $Y$ lies in the essential image. We start with showing this for $Y^{(1)}$.
We know that there are Eilenberg-MacLane objects $E_{I^1}, E_{I^2} \in \Da$ such that $\R(E_{I^1}) \cong \mathcal{E}_{I^1}$
and $\R(E_{I^2}) \cong \mathcal{E}_{I^2}$.
We started with an injective resolution
\[
E(1)_*(Y) \rightarrow I^0 \rightarrow I^1 \rightarrow I^2 \rightarrow 0
\]
for $E(1)_*(Y) \in \A$. Using Lemma \ref{hocolim}, this resolution equals
\begin{equation}
\label{res}
E(1)_*(Y) \rightarrow \bigoplus\limits_i H^i(E_{I^0})[-i] \stackrel{d^1}{\rightarrow} \bigoplus\limits_i H^i(E_{I^1})[-i]
\stackrel{d^2}{\rightarrow} \bigoplus\limits_i H^i(E_{I^2})[-i] \rightarrow 0
\end{equation}
with above Eilenberg-MacLane objects in $\Da$. We take those Eilenberg-MacLane objects and complete them to an exact triangle
\begin{equation}
E_{I^2} \rightarrow D \rightarrow E_{I^1} \rightarrow E_{I^2}[1]
\end{equation}
in $\Da$. Applying
\[
\bigoplus\limits_i H^i(-)[-i]
\]
 to this triangle gives a long exact sequence in $\A$. Since $d^2$ in
(\ref{res}) is a surjection, the third morphism in this triangle
induces a surjection in cohomology as well. Consequently, the second
morphism $D \rightarrow E_{I^1}$ must give an injection in
cohomology. So we can apply the formerly stated fact at the end of
section 2 again that
\[
\R(E_{I^2}) \rightarrow \R(D) \rightarrow \R(E_{I^1}) \rightarrow \R(E_{I^2}[1])
\]
is an exact triangle in $\Ho(\E)$.

Consider

\[
\xymatrix{ \mathcal{E}_{I^2} \ar[r]\ar[d]_{\cong}^{\R}  & Y^{(1)} \ar[r]\ar@{.>}[d] & \mathcal{E}_{I^1} \ar[r]\ar[d]_{\cong}^{\R} & \Sigma\mathcal{E}_{I^2} \ar[d]_{\cong}^{\R} \\
\R(E_{I^2}) \ar[r] & \R(D^*) \ar[r] & \R(E_{I^1}) \ar[r] & \R(E_{I^2}[1])
}
\]

with the upper triangle coming from (\ref{Yadamsres}).
The third square commutes since $\R$ is full. By the axioms of a triangulated category there exist a morphism
$Y^{(1)} \rightarrow \R(D^*)$ making the whole diagram commute. By the 5-lemma, this is an isomorphism, thus $Y^{(1)} \cong \R(D^*)$, and so $Y^{(1)}$ lies in the essential image of $\R$.
Similarly, this also follows for $Y$, which completes the proof that $\R$ is an equivalence of categories.
\end{proof}

\begin{ccc}
$\R$ preserves the Adams filtration.
\end{ccc}

\begin{rrr}
Nora Ganter recently proved in  that for the case of $E(1)$-local spectra $\R$ is not just an equivalence of categories
but $\R$ also carries tensor products of cochain complexes into smash
products of spectra. (This is not known to be true for arbitrary $n$ with $n^2+n < 2p-2$.)
\end{rrr}

\section{A further application}\label{quillen}

As proved, $\R$ provides an abstract equivalence of triangulated categories which also happen to be
homotopy categories of model categories. The next question now is if $\Db$ and
$\Ho(\E)$ are equivalent as categories, can their model structures also be positively compared,
i.e. is there a Quillen equivalence between them?

The answer to that is remarkable:

\begin{ppp}
The categories $\Db$ and $\Ho(\E)$ are not Quillen equivalent. In particular, $\R$ is not a Quillen equivalence.
\end{ppp}

\begin{proof}
To prove this, we compare the homotopy types of certain mapping spaces for each category. Let us first collect the necessary definitions.
For a pointed simplicial model category $\C$ there is a mapping space functor
\[
\map_{\C}(-,-): \C^{op} \times \C \longrightarrow \sset
\]

to the category of pointed simplicial sets satisfying
\[
\map_{\C}(X,Y)_0 = \Hom_{\C}(X,Y)
\]
for all $X, Y \in \C$ and certain adjointness properties (see e.g. \cite{GJ99}, Definition II.2.1).
However, $\Da$ and $\Db$ are not simplicial categories. The next best thing we can achieve is a notion of a mapping space that is well-defined up to homotopy, which will do for our purposes.

To achieve this, we look at the category $\C^{\Delta}$ of cosimplicial objects in $\C$ and view $X$ as constant object in $\C^{\Delta}$. The category $\C^{\Delta}$ of cosimplicial objects in a model category $\C$ can be given a model structure, the so-called Reedy model structure. For details of this, see \cite{Hov99} Section 5.2. We now define a special replacement of $X$ in $\C^{\Delta}$, so-called {\it frames}. To do this, we first need the following:

\begin{ddd}
Via the methods of \cite{Hov99}, Remark 5.2.3. and Example 5.2.4., there are functors $\fatl^{\bullet}, \fatr^{\bullet}: \C \longrightarrow \C^{\Delta} $ with the following properties:
\newline
Let $X \in \C$:
\begin{itemize}
\item the $n^{th}$ level space of the object $\fatl^{\bullet}X$ is the $n+1$-fold coproduct of $A$
\item $\fatl^{\bullet}: \C \longrightarrow \C^{\Delta}$ is a left adjoint to the evaluation functor $ev_0: \C^{\Delta} \longrightarrow \C$ that sends $A^{\bullet}$ to $A^{\bullet}[0]$
\item the $n^{th}$ level space of the object $\fatr^{\bullet}X$ is $X$ itself
\item $\fatr^{\bullet}: \C \longrightarrow \C^{\Delta}$ is a right adjoint to $ev_0: \C^{\Delta} \longrightarrow \C$
\end{itemize}
\end{ddd}

\begin{rrr} One can prove that $\fatr^{\bullet}$ is the constant cosimplicial functor.
There is a natural transformation $\fatl^{\bullet} \longrightarrow \fatr^{\bullet}$ that is the identity in degree zero and the fold map in higher degrees.
\end{rrr}

With these functors, we can now define cosimplicial frames:
\begin{ddd}
Let $\C$ be a model category, $X$ an object of $\C$. A {\it cosimplicial frame} for $X$ is a cosimplicial object $X^{\bullet} \in \C^{\Delta}$ together with a factorisation of the map $\fatl^{\bullet}X \longrightarrow \fatr^{\bullet}X$ in $\C^{\Delta}$
$$\xymatrix@M+2pt{
\fatl^{\bullet}X \ar@{>->}[r] & X^{\bullet} \ar[r]^{\sim} & \fatr^{\bullet}X}
$$
where the weak equivalence $X^{\bullet} \stackrel{\sim}{\longrightarrow} \fatr^{\bullet}X$ in degree zero induces a weak equivalence in $\C$.
\end{ddd}

For the existence of such framings, see \cite{Hov99}, Theorem 5.2.8.

We now use this definition to define mapping spaces:

\begin{ddd}\label{mapdef}
Let $X, Y$ be objects of $\C$, $X^{\bullet}$ a cosimplicial frame for $X$ and
\[
\xymatrix@M+2pt{ Y \ar@{>->}[r]^{\sim} & Y^{\fib} \ar@{>>}[r] & \ast
}
\]
a factorisation of $Y \rightarrow \ast$.
Then the {\it (left) mapping space} for $X$ and $Y$ is defined via
\[
\map_{\C}(X,Y):= \C(X^{\bullet},Y^{\fib}) \in \sset,
\]
where $\C(X^{\bullet},Y^{\fib})$ is the simplicial set with
\[
\C(X^{\bullet},Y^{\fib})_n := \Hom_{\C}(X^{\bullet}[n], Y^{\fib}).
\]
\end{ddd}

However, it is not clear whether this definition actually deserves to be called a definition since it depends on two choices: firstly, the cosimplicial frame for $X$ and secondly, the fibrant replacement for $Y$. So, for this definition to make sense we need the following:

\begin{lem}\label{framingchoice}
Let $X_1^{\bullet}, X_2^{\bullet}$ be two cosimplicial frames for cofibrant $X$ in $\C$, and let $Y_1^{\fib}, Y_2^{\fib}$ be two fibrant replacements for $Y$. Then
\[
\C(X_1^{\bullet}, Y_1^{\fib}) \simeq \C(X_2^{\bullet}, Y_2^{\fib})
\]
in $\sset$.
\end{lem}

\begin{proof}
First, let $X_1^{\bullet}$ and $X_2^{\bullet}$ be two cosimplicial frames for $X$. By definition, the frames $X_1^{\bullet}$ and $X_2^{\bullet}$ are linked by a zig-zag of weak equivalences
\[
X_1^{\bullet} \stackrel{\sim}{\longrightarrow} \fatr^{\bullet}X \stackrel{\sim}{\longleftarrow} X_2^{\bullet}.
\]
For fibrant $Y$, the functor $\C(-,Y)$ preserves weak equivalences (\cite{SchShi02} Lemma 6.3), so for fibrant $Y$ and $X_1^{\bullet}$, $X_2^{\bullet}$ as above, we have
\[
\C(X_1^{\bullet}, Y) \simeq \C(X_2^{\bullet}, Y).
\]

For the second part we quote \cite{Hov99}, Corollary 5.4.4, which says that for fibrant $X$ in $\C$, the functor
\[
\C(X^{\bullet},-): \C \longrightarrow \sset
\]
preserves fibrations and acyclic fibrations, in particular between fibrant objects. So Ken Brown's lemma applies
 (see e.g. \cite{Hov99}, Lemma 1.1.12), and it follows that $\C(X^{\bullet},-)$ takes weak equivalences between fibrant objects in $\C$ to weak equivalences in $\sset$ which proves the claim of our lemma.
\end{proof}

Now we look at the behaviour of mapping spaces under Quillen functors and Quillen equivalences.

\begin{lem}\label{mapQuillen}
Let $L: \C \rightleftarrows \D : R$ be a Quillen equivalence, $X, X' \in \C$ both cofibrant. Then
\[
\map_{\C}(X,X') \cong \map_{\D}(LX,LX')
\]
in Ho($\sset$).
\end{lem}

\begin{proof}
First of all, let $L: \C \rightleftarrows \D : R$ be a Quillen adjoint functor pair, $X \in \C$ and $Y \in \D$. Then
\[
\map_{\D}(LX,Y) = \D((LX)^{\bullet}, Y^{\fib})
\]
by definition. Since $L$ is a left Quillen functor, $L(X^{\bullet}) \in \D^{\Delta}$ is also a cosimplicial frame for $LX$ (\cite{Hov99}, Lemma 5.6.1),
so
\[
\D((LX)^{\bullet}, Y^{\fib}) \cong \D(L(X^{\bullet}), Y^{\fib})
\]
by Lemma \ref{framingchoice}. By adjointness,
\[
\Hom_{\D}(L(X^{\bullet})[n], Y^{\fib}) \cong \Hom_{\C}(X^{\bullet}[n],R(Y^{\fib})),
\]
so
\[
\D(L(X^{\bullet}), Y^{\fib}) \cong \C(X^{\bullet}, R(Y^{\fib})).
\]
Since $R$ is a right Quillen functor, $R(Y^{\fib})$ is a fibrant replacement for $RY$, consequently by Lemma \ref{framingchoice},
\[
\C(X^{\bullet}, R(Y^{\fib})) \simeq \C(X^{\bullet}, (RY)^{\fib}) = \map_{\C}(X,RY).
\]
Thus, altogether we have
\begin{equation}
\label{mapad}
\map_{\C}(X,RY) \simeq \map_{\D}(LX,Y).
\end{equation}

Next, let $L:\C \rightleftarrows \D:R$ be a Quillen equivalence and $X' \in \C$ cofibrant. Then
\[
LX' \stackrel{\sim}{\longrightarrow} (LX')^{\fib}
\]
is a weak equivalence in $\D$ with cofibrant source and fibrant target, so by definition of a Quillen equivalence, the adjoint map
\[
X' \stackrel{\sim}{\longrightarrow} R((LX')^{\fib})
\]
is a weak equivalence in $\C$. Since $R$ is a right Quillen functor, $R((LX')^{\fib})$ is fibrant in $\C$. Consequently, $R((LX')^{\fib})$ is a fibrant replacement for $X'$ in $\C$. By Lemma \ref{framingchoice} and above adjointness result for mapping spaces (\ref{mapad}), it follows that
\[
\map_{\C}(X,X') \simeq \map_{\C}(X,R((LX')^{\fib})) \simeq \map_{\D}(LX,LX')
\]
in $\sset$ which proves the lemma.
\end{proof}

Back to our special case: We will see that for all $C, D \in \Cb$, $\map_{\Cb}(C,D)$ is weakly equivalent to a product of Eilenberg-MacLane spaces. However, the mapping space $\map_{\E}(S^0,S^0)$ is not a product of Eilenberg-MacLane spaces, so as a consequence of Lemma \ref{mapQuillen}, there is no Quillen equivalence between those two model categories which was the claim of the proposition.

The category $\Cb$ is abelian, so for all $C_1, C_2 \in \Cb$,  the $n$-simplices of $\map_{\Cb}(C_1,C_2)$
\[
\C(C_1^{\bullet},C_2^{\fib})_n = \Hom(C_1^{\bullet}[n],C_2)
\]
form an abelian group, and the simplicial structure maps are group homomorphisms, so
\[
\C(C_1^{\bullet},C_2^{\fib}) = \map_{\Cb}(C_1,C_2)
\]
is not just a simplicial set but a simplicial abelian group. From Proposition III.2.20 of \cite{GJ99}, it follows that
\[
\map_{\Cb}(C_1,C_2) \cong \prod\limits_{n \ge 0} K(\pi_n \map_{\Cb}(C_1,C_2)_n , n)
\]
where $K(G,n)$ denotes the $n^{th}$ Eilenberg-MacLane space for the abelian group $G$.

However, there are spectra for which the mapping spaces over $\E$ are not products of Eilenberg-MacLane spaces, for example $\map_{\E} (S^0,S^0) \cong QL_1 S^0 = \colim_n \Omega^n L_1 S^n.$ Thus, $\Cb$ and $\E$ cannot be Quillen equivalent and $\Cb$ provides an exotic model for $\E$.

\end{proof}

In other words, $\Cb$ provides an exotic model for $\E$. For the stable homotopy
category itself such exotic models do not exist, as proved by Schwede in \cite{Sch05}. However, this is not true for the chromatic localisations of the stable homotopy
category in the cases $n^2+n < 2p-2$ (shown here explicitly for $n=1$). It is not yet known how many such exotic models exist
and what can be said about the other chromatic localisations.


\begin{thebibliography}{Hov99}

\bibitem[BF78]{BouFri78}
A.~K. Bousfield and E.~M. Friedlander.
\newblock Homotopy theory of {$\Gamma $}-spaces, spectra, and bisimplicial
  sets.
\newblock In {\em Geometric applications of homotopy theory (Proc. Conf.,
  Evanston, Ill., 1977), II}, volume 658 of {\em Lecture Notes in Math.}, pages
  80--130. Springer, Berlin, 1978.

\bibitem[Bou85]{Bou85}
A.~K. Bousfield.
\newblock On the homotopy theory of {$K$}-local spectra at an odd prime.
\newblock {\em Amer. J. Math.}, 107(4):895--932, 1985.

\bibitem[Fra96]{Fra96}
J.~Franke.
\newblock Uniqueness theorems for certain triangulated categories possessing an
  adams spectral sequence.
\newblock {\em http://www.math.uiuc.edu/K-theory/0139/}, 1996.

\bibitem[GJ99]{GJ99}
P.~G. Goerss and J.~F. Jardine.
\newblock {\em Simplicial homotopy theory}, volume 174 of {\em Progress in
  Mathematics}.
\newblock Birkh\"auser Verlag, Basel, 1999.

\bibitem[Hov99]{Hov99}
M.~Hovey.
\newblock {\em Model categories}, volume~63 of {\em Mathematical Surveys and
  Monographs}.
\newblock American Mathematical Society, Providence, RI, 1999.

\bibitem[Sch05]{Sch05}
S.~Schwede.
\newblock The stable homotopy category is rigid.
\newblock {\em http://www.math.uni-bonn.de/people/schwede/rigid.pdf}, 2005.

\bibitem[SS02]{SchShi02}
S.~Schwede and B.~Shipley.
\newblock A uniqueness theorem for stable homotopy theory.
\newblock {\em Math. Z.}, 239(4):803--828, 2002.

\end{thebibliography}
\end{document}